\documentclass[11pt]{amsart}
\begin{document}
\newtheorem{thm}{Theorem}[section]
\newtheorem{prop}[thm]{Proposition}
\newtheorem{lem}[thm]{Lemma}
\newtheorem{cor}[thm]{Corollary}
\newtheorem{question}[thm]{Question}

\theoremstyle{definition}
\newtheorem{defn}[thm]{Definition}
\newtheorem{remark}[thm]{Remark}
\newtheorem{example}[thm]{Example}
\newtheorem{reduction}[thm]{Reduction}

\newcommand{\Span}{\operatorname{Span}}
\newcommand{\hr}{\operatorname{hr}} 
\newcommand{\Gr}{\operatorname{Gr}}
\newcommand{\Id}{\operatorname{Id}}
\newcommand{\id}{\operatorname{id}}
\newcommand{\pr}{\operatorname{pr}}
\newcommand{\Hom}{\operatorname{Hom}}
\newcommand{\hw}{\operatorname{hw}}

\newcommand{\bbA}{{\mathbb A}}
\newcommand{\bbC}{{\mathbb C}}
\newcommand{\bbR}{{\mathbb R}}
\newcommand{\bbZ}{{\mathbb Z}}
\newcommand{\bbP}{{\mathbb P}}
\newcommand{\bbQ}{{\mathbb Q}}

\newcommand{\kbar}{\overline k}

\newcommand{\mo}{\mathopen<}
\newcommand{\mc}{\mathclose>}

\numberwithin{equation}{section}

\newcommand{\Ker}{\operatorname{Ker}}
\newcommand{\Aut}{\operatorname{Aut}}
\newcommand{\sdp}{\mathbin{{>}\!{\triangleleft}}} 
\newcommand{\Alt}{\operatorname{A}}   
\newcommand{\GL}{\operatorname{GL}}
\newcommand{\PGL}{\operatorname{PGL}}
\newcommand{\SL}{\operatorname{\rm SL}}
\newcommand{\rank}{\operatorname{rank}}
\newcommand{\aut}{\operatorname{Aut}}
\newcommand{\Char}{\operatorname{\rm char\,}} 
\newcommand{\diag}{\operatorname{\rm diag}}
\newcommand{\lra}{\longrightarrow}
\newcommand{\SO}{\operatorname{SO}}
\newcommand{\M}{\operatorname{M}}        
\newcommand{\cO}{\mathcal{O}}
\newcommand{\cA}{\mathcal{A}}
\newcommand{\cB}{\mathcal{B}}
\newcommand{\cC}{\mathcal{C}}
\newcommand{\cE}{\mathcal{E}}
\newcommand{\cT}{\mathcal{T}}
\newcommand{\Sym}{{\operatorname{S}}}    
\newcommand{\tr}{\operatorname{\rm tr}}
\newcommand{\trace}{\tr}

\title[Nesting maps]{Nesting maps of Grassmannians} 
\author{C. De Concini}
\address{Dipartimento di Matematica, Universit\'a di Roma ``La Sapienza",
Piazzale Aldo Moro 5, 00185 ROMA, Italy}
\email{deconcin@mercurio.mat.uniroma1.it}

\author{Z. Reichstein}
\address{Department of Mathematics, University of British Columbia,
        Vancouver, BC V6T 1Z2, Canada}
\email{reichste@math.ubc.ca}
\urladdr{www.math.ubc.ca/$\stackrel{\sim}{\phantom{.}}$reichst}
\thanks{Z. Reichstein was supported in part by an NSERC research grant}


%


\begin{abstract} 
Let $F$ be a field and $\Gr(i, F^n)$ be the Grassmannian of
$i$-dimensional linear subspaces of $F^n$. A map 
$f \colon \Gr(i, F^n) \lra \Gr(j, F^n)$ is called nesting
if $l \subset f(l)$ for every $l \in \Gr(i, F^n)$. Glover, 
Homer and Stong showed that there are no continuous 
nesting maps $\Gr(i, \bbC^n) \lra \Gr(j, \bbC^n)$ 
except for a few obvious ones. We prove a similar result
for algebraic nesting maps $\Gr(i, F^n) \lra \Gr(j, F^n)$, 
where $F$ is an algebraically closed field of arbitrary 
characteristic. For $i = 1$ this yields a description 
of the algebraic subbundles of the tangent bundle 
to the projective space $\bbP_{F}^n$. 
\end{abstract}

\maketitle

\section{Introduction}
\label{sect.intro}
Let $F$ be a field. We shall denote the Grassmannian of $i$-dimensional 
linear subspaces of $F^n$ by $\Gr(i, F^n)$. Suppose
$i$, $j$ and $n$ are integers satisfying 
$1 \leq i < j \leq n-1$.  We shall say that a map of Grassmannians
\begin{equation} \label{e1}
f \colon \Gr(i, F^n) \lra \Gr(j, F^n)
\end{equation}
is {\em nesting}, if $l \subset f(l)$ for every $l \in \Gr(i, F^n)$.
The starting point for this note is the following combinatorial
result communicated to us by J. Buhler.

\begin{thm} \label{thm0}
Let $F$ be a finite field and
$i < j$ be positive integers such that $n = i + j$.
Then there exists a bijective nesting map (of sets)
$f \colon \Gr(i, F^n) \lra \Gr(j, F^n)$.
\end{thm}

Theorem~\ref{thm0} can be deduced from a theorem 
of P. Hall about systems of distinct representatives; 
for details see Section~\ref{sect2}.

It is natural to ask for what $n$, $i$ and $j$
there exist {\em algebraic} (and for $F = \bbR$, 
$F = \bbC$, {\em continuous}) nesting 
maps $f \colon \Gr(i, F^n) \lra \Gr(j, F^n)$.

\begin{example} \label{ex1}
Let $n$ be an even integer and $\omega$ be an alternating
bilinear form on $F^n$, i.e., a non-degenerate 
bilinear form $\omega(v, w)$ such that $\omega(v, v) = 0$ for 
every $v \in F^n$. (If $\Char(F) \neq 2$, then the last condition 
is equivalent to $\omega(v, w) = - \omega(w, v)$ for every 
$v, w \in F^n$, so ``alternating" is the same as ``symplectic".) 
Then $f \colon \Gr(1, F^n) \lra \Gr(n-1, F^n)$, given by 
$f(l) = l^{\perp_{\omega}}$, is an algebraic nesting isomorphism. 
Here $l^{\perp_{\omega}}$ is the orthogonal complement with 
respect to $\omega$.
\qed
\end{example}

\begin{example} \label{ex4.1} Suppose $n \geq 2$ is even,
$\omega$ is an alternating form and $H$ is
a positive-definite Hermitian form on $\bbC^n$. Then
we can define a {\em continuous} nesting map
$g \colon \Gr(1, \bbC^n) = \bbP_{\bbC}^{n-1} \lra \Gr(2, \bbC^n)$
by  $l \mapsto l \oplus (l^{\perp_{\omega}})^{\perp_H}$.

Note that if $n$ is odd then there is no continuous map
$g \colon \bbP_{\bbC}^{n-1} \lra \Gr(2, \bbC^n)$. This is a special 
case of Theorem~\ref{thm1b}(a) below, but it can also be seen directly 
as follows. Suppose
$g$ exists. Then $\alpha \colon l \mapsto l^{\perp_H} \cap g(l)$
is a continuous fixed point free map $\bbP_{\bbC}^{n-1} \lra
\bbP_{\bbC}^{n-1}$, contradicting the Lefschetz fixed 
point theorem; see,~e.g.,~\cite[30.13]{greenberg}. 
(Note that $\alpha$ is, indeed, well-defined, 
because $l \subset g(l)$ implies $g(l) \not \subset l^{\perp_H}$.)
\qed
\end{example}

The following results assert that, over an algebraically closed
field $F$, these are essentially the only examples.
First consider the case where $i \geq 2$. 

\begin{thm} \label{thm1a} Suppose $2 \leq i < j \leq n-1$. Then

\smallskip
(a) there does not exist a continuous nesting
map $\Gr(i, \bbC^n) \lra \Gr(j, \bbC^n)$,

\smallskip
(b) there does not exist an algebraic nesting
map $\Gr(i, F^n) \lra \Gr(j, F^n)$ for any
algebraically closed field $F$.
\end{thm}

For $i = 1$ the continuous and the algebraic cases diverge. We shall
write $\bbP_{F}^{n-1}$ in place of $\Gr(1, F^n)$.

\begin{thm} \label{thm1b} Suppose $2 \leq j \leq n-1$. Then

\smallskip
(a) a continuous nesting map 
$\bbP_{\bbC}^{n-1} \lra \Gr(j, \bbC^n)$ exists
if and only if $n$ is even and $j = 2$ or $j = n-1$.

\smallskip
(b) Let $F$ be an algebraically closed field and
$f \colon \bbP_F^{n-1} \lra \Gr(j, F^n)$ be an algebraic 
nesting map. Then $n$ is even, $j = n-1$, and $f$ is 
as in Example~\ref{ex1}.
\end{thm}

Nesting maps $\bbP_{F}^{n-1} \lra \Gr(j, F^n)$ are easily
seen to be in a natural 1-1-correspondence with rank $j-1$
subbundles of $T(\bbP_F^{n-1})$; cf. Lemma~\ref{lem4.1}.
Using this correspondence, we can rephrase 
Theorem~\ref{thm1b} as follows:

\begin{cor} \label{cor2} Let $F$ be an algebraically closed field and
$T(\bbP_{F}^{n-1})$ be the tangent bundle to the projective
space $\bbP_{F}^{n-1}$.  

\smallskip
(a) $T(\bbP_{\bbC}^{n-1})$ has a topological (complex) 
subbundle of rank $1 \leq r \leq n-2$ if and only if $n$ 
is even and $r = 1$ or $n-2$. 

\smallskip
(b) If $n$ is odd then $T(\bbP_{F}^{n-1})$ has no nontrivial algebraic 
subbundles.  If $n$ is even then the only non-trivial algebraic subbundle of
$T(\bbP_{F}^{n-1})$ (up to an automorphism of $\bbP_F^{n-1}$) 
is the null correlation bundle of rank $n - 2$.
\end{cor}

Here the null-correlation bundle on $\bbP_{F}^{n-1}$ (defined,
e.g., in~\cite[Section 4.2]{oss} or~\cite[p. 128]{bott}) 
is the rank $n-2$ bundle associated to the nesting map of Example~\ref{ex1}.  
(Note that different choices of the alternating form in Example~\ref{ex1}
give rise to bundles that differ by an automorphism of $\bbP_{F}^{n-1}$.)

Theorem~\ref{thm1a}(a) is due to Stong~\cite[Theorem 2]{stong}.
For a field $F$ of characteristic zero, part (b) follows 
from part (a) by the Lefschetz priniciple. 
Theorem~\ref{thm1b}(a) and Corollary~\ref{cor2}(a) are due to 
Glover, Homer and Stong~\cite[Theorem 2(ii) and Theorem 1.1(ii)]{ghs}. 
(A special case of Corollary~\ref{cor2}(a) was considered earlier 
by Bott~\cite[Corollary 1.5]{bott}.)
For fields $F$ of characteristic zero, part (b) can be deduced
from these results and the work of Roan~\cite[Main Theorem]{roan}.

The purpose of this note is to give uniform characteristic-free 
proofs of Theorem~\ref{thm1a}, Theorem~\ref{thm1b} 
and Corollary~\ref{cor2}. We also include a proof 
of Theorem~\ref{thm0} in a short appendix. 
An application of Theorem~\ref{thm1b} can be found in a recent
paper of Bergman~\cite{bergman}.

%

\section*{Acknowledgments} We are grateful to J. Buhler for
contributing Theorem~\ref{thm0} and to G. Bergman, 
J. Carrell, M.-C. Chang, K.-Y. Lam, S. L$'$vovsky, 
M. Thaddeus and A. Vistoli for stimulating discussions.

\section{Proof of Theorem~\ref{thm1a}}
\label{sect.i>=2}

We begin by recalling some well-known facts about
the cohomology of the Grassmannian; for details we 
refer the reader to~\cite[Chapter 14]{fulton-it} (see also~\cite{carrell}
or \cite{hiller}).

Let $\sigma_i$ be the $i$th elementary symmetric polynomial
in the independent variables $x_1, \dots, x_n$.
The cohomology ring $H^*(\Gr(i, F^n), \bbZ)$
is isomorphic to the quotient of $\bbZ[\sigma_1,\ldots \sigma_i]$, by the
ideal generated by the Schur functions $S_{\lambda}$, where
$\lambda$ ranges over the partitions of $n$ whose Young tableaux are
not contained in a rectangle
with $i$ rows and $n-i$ columns. Here, by convention,
$\sigma_i$ is the Schur function for a column with $i$ rows.

Two remarks are now in order. First of all, the cohomology class
corresponding to $\sigma_i$ lies in $H^{2i}(\Gr(i, F^n), \bbZ)$;
however, for notational convenience, we will continue
to use the natural grading induced from $\bbZ[x_1, \dots, x_n]$, so that
$\deg (\sigma_i) = i$. Secondly, since $\bbZ[\sigma_1,\ldots \sigma_i]$
is generated over $\bbZ$ by the Schur functions $S_{\lambda}$,
where $\lambda$ ranges over partitions with at most $i$ rows,
we identify $H^*(\Gr(i, F^n), \bbZ)$ with $\bbZ[\sigma_1,\ldots \sigma_i]/I$, 
where $I$ is a homogeneous ideal generated by elements of degree $> n-i$. 

We will prove both parts of Theorem~\ref{thm1a} by the same 
computation. (If $\Char(F) >0$ in part (b), we replace 
the cohomology ring $H^*(\Gr(i, F^n), \bbZ)$ by the Chow ring;
our computation will remain valid there; cf.~\cite[Chapter 14]{fulton-it}.)

\smallskip
Let $V$ be the trivial bundle of rank $n$ on $\Gr(i, F^n)$ and let
$\cT_i$ be the tautological bundle of rank $i$. (The fiber 
of $\cT_i$ over $l \in \Gr(i, F^n)$ consists of the vectors in $l$.)
Since $f$ is a nesting map,
\[ \text{$\cT_i \subset f^{\ast}(\cT_j)$ and
$f^{\ast}(\cT_j) \subset V$,} \]
where $\subset$ means ``subbundle of". In other words, 
\begin{equation} \label{e.AB}
\text{$\cB = f^{\ast}(\cT_j)/\cT_i$ is a subbundle of
$\cA = V/\cT_i$,}
\end{equation}
where $\rank(\cA) = n - i$ and $\rank(\cB) = j - i$.

Now recall that the $r$th Chern class of the tautological bundle 
is given by $c_r(\cT_i) = (-1)^r \sigma_r$ (cf., e.g.,~\cite[p. 111]{hiller}),
so that the total Chern class is 
\[ c_{\rm tot}(\cT_i) = 
(1 - \sigma_1 + \sigma_2 - \dots + (-1)^i \sigma_i) \, . \] 
Since  $\cA:= V/\cT_i$, we have
\begin{equation} \label{e.chern}
c_{\rm tot}(\cA) (1 - \sigma_1 + \sigma_2 - \dots + (-1)^i \sigma_i) =  1 \, ,
\end{equation}
where $c_{\rm tot}(\cA)$ is the total Chern class of $\cA$, and
$1$ represents the total Chern class of the trivial bundle $V$
on $\Gr(i, F^n)$. 
 
Since $\cA$ has a subbundle $\cB$ or rank $n -j$, 
$c_{\rm tot}(\cA)$ factors as a product of elements
of degree $n-j$ and $j - i$ in $H^{\ast}(\Gr(i, F^n), \bbZ)$.
We will show that this is impossible. Consider the
degree-preserving homomorphism
\[ \phi \colon \bbZ[x_1, \dots, x_i, \dots, x_n] \lra \bbZ[x_1, \dots, x_i] \]
given by $\phi(x_r) = x_r$ if $1 \leq r \leq i$ and
$\phi(x_r) = 0$ if $i + 1 \leq r \leq n$. Under this homomorphism
$\bbZ[\sigma_1, \dots, \sigma_i]$ maps isomorphically to
$\bbZ[\phi(\sigma_1), \dots, \phi(\sigma_i)]$,
where $\phi(\sigma_r)$ is the $r$th elementary symmetric 
function in $x_1, \dots, x_i$. Thus it is enough to prove that
$\phi(c_{\rm tot}(\cA))$ is irreducible in 
$\bbZ[\phi(\sigma_1), \dots, \phi(\sigma_i)]$. We will, in fact, prove that
$\phi(c_{\rm tot}(\cA))$ is irreducible in $\bbZ[x_1, \dots, x_i]$ and
even in $\bbC[x_1, \dots, x_i]$.
Applying $\phi$ to both sides of~\eqref{e.chern}, we see that
\begin{gather} 
1 = \phi(c_{\rm tot}(\cA)) 
(1 - \phi(\sigma_1) + \phi(\sigma_2) - \dots + (-1)^i \phi(\sigma_i)) = 
\nonumber \\
= \phi(c_{\rm tot}(\cA)) (1 - x_1)(1-x_2) \dots (1-x_i) \nonumber \, . 
\end{gather}
in $\bbZ[x_1, \dots, x_i]/\phi(I)$. Since $\phi(c_{\rm tot}(\cA))$ has
degree $\leq n - i$, and $I$ (and thus $\phi(I)$) is generated 
by homogeneous elements of degree $> n - i$, 
\begin{gather}
\phi(c_{\rm tot}(\cA)) = \bigl( (1 - x_1)^{-1} (1-x_2)^{-1} \dots 
(1-x_i)^{-1} \bigr)_{|n-i} =  \label{e2.5} \\ 
\bigl( (1 + x_1 + x_1^2 + \dots)(1 + x_2 + x_2^2 + \dots) \dots
(1 + x_i + x_i^2 + \dots) \bigr)_{|n-i} = \nonumber \\
1 + \theta_1(x_1, \dots, x_i) +
\theta_2(x_1, \dots, x_i) + \dots + \theta_{n-i}(x_1, \dots, x_i) \, , 
\nonumber
\end{gather}
where ${ }_{|n-i}$ means that we cut the series at the terms of degree at most
$n-i$ and  $\theta_d(x_1, \dots, x_i)$ 
denotes the sum of all monomials of degree $d$ in $x_1, \dots, x_i$.  
It remains to show that
$1 + \theta_1(x_1, \dots, x_i) + \dots + \theta_{n-i}(x_1, \dots, x_i)$
is irreducible in $\bbC[x_1,\ldots, x_i]$. Homogenizing
this polynomial with respect to an additional variable $x_0$, we obtain
$\theta_{n-i}(x_0,x_1,\dots,x_i)$. Thus Theorem~\ref{thm1a}
is a consequence of the following:

\begin{lem} \label{lem.s_d} $\theta_{d}(x_0,x_1,\ldots ,x_i)$ is
irreducible in $\bbC[x_0, \dots, x_i]$ for any $i \geq 2$ and $d\geq 1$.
\end{lem}

\begin{proof} It suffices to consider the case $i = 2$.
For notational convenience, we will write $x, y$ and $z$ for $x_0$, $x_1$
and $x_2$. In other words, we want to prove that
the projective curve $X_d \subset \bbP^2$ given by
$\theta_{d}(x,y,z)=0$ is irreducible. It suffices to show that
$X_d$ is non-singular.  
We proceed by induction on $d$.

The base case, $d=1$, is trivial, since $\theta_{1}(x,y,z)=x+y+z$ cuts out
a line in $\bbP^2$. The inductive step will rely on the formulas
\begin{equation} \label{e3.2}
\theta_{d}(x,y,z)=x\theta_{d-1}(x,y,z)+\theta_{d}(y,z)
\end{equation}

\noindent
and
\begin{equation} \label{e3.3}
(\frac{\partial }{\partial x} +  \frac{\partial}{\partial y} +
\frac{\partial}{\partial z}) \theta_{d}(x,y,z)=(d+2) \theta_{d-1}(x,y,z) \, .
\end{equation}
Suppose $p = (x_0: y_0: z_0)$ is a singular point of $X_d$.  By~\eqref{e3.3},
\begin{equation} \label{e3.4}
\text{$\theta_{d-1}(p)=0$, i.e., $p \in X_{d-1}$.}
\end{equation} 
We claim that
$x_0 y_0 z_0 \neq 0$. Indeed, assume the contrary, say, $p = (0: 1: z_0)$.
Then formulas~\eqref{e3.2}, and~\eqref{e3.4} say that
$\theta_{d}(1,z_0)=0$, i.e.,
$z_0 \neq 1$ is a $(d+1)$th root of unity. On the other hand,
by~\eqref{e3.4},
\[ 0 = \theta_{d-1}(p) = \theta_{d-1}(0, 1, z_0) = \theta_{d-1}(1, z_0) \, , \]
so that $z_0 \neq 1$ is a $d$th root of unity, a contradiction. This
proves the claim.

Differentiating~\eqref{e3.2} with respect to $x$, we obtain
$$\frac{\partial}{\partial x} \theta_{d}(x,y,z)= \theta_{d-1}(x,y,z)+
x\frac{\partial}{\partial x} \theta_{d-1}(x,y,z) \, ,$$
and similarly for $y$ and $z$.
Combining this with~\eqref{e3.2} and~\eqref{e3.4}, and keeping in mind that
$x_0 y_0 z_0 \neq 0$, we obtain
$$ \theta_{d-1}(p)= \frac{\partial}{\partial x} \theta_{d-1}(p)=
\frac{\partial}{\partial y} \theta_{d-1}(p)=
\frac{\partial}{\partial z} \theta_{d-1}(p)=0 \, . $$
Thus $p$ is a singular point of $X_{d-1}$, contradicting
our induction assumption.
This completes the proof of Lemma~\ref{lem.s_d} and of Theorem~\ref{thm1a}.
\end{proof}

\section{The Schwarzenberger conditions}

Our proof of Theorem~\ref{thm1b} in the next section will rely 
on the following result about vector bundles on projective spaces.

\begin{prop} \label{prop.schw}
Suppose $E$ is a continuous vector bundle of rank $1 \leq r \leq m-2$ 
on $\bbP_{\bbC}^m$ or an algebraic vector bundle
on $\bbP_{F}^m$, where $F$ is an algebraically closed field.
Let $p(t) = 1 + c_1 t + \dots + c_r t^r$ be the Chern polynomial of $E$.

Suppose the (complex) roots $w_1, \dots, w_r$ of $p(t)$ 
are distinct and lie on the unit circle.  Then either $r = 1$ and 
$p(t) = 1 \pm t$ or $r = 2$ and $p(t) = 1 - t^2$.
\end{prop}

Note that by a theorem of Kronecker, $w_1, \dots, w_r$ are 
necessarily roots of unity; however, we shall not use this  
in the proof.

\begin{proof}
Since $w_1, \dots, w_r$ lie on the unit circle 
each $w_i^{-1} = \overline{w_i}$ is also a root of $p(t)$, so that
\[ p(t) = (1 - w_1t) \dots (1 - w_r t) \, . \]

Let $B_{k, l} = \sum_{i= 1}^r \begin{pmatrix} k - w_i \\ l 
\end{pmatrix}$. 
Our argument is based on the Schwarzenberger conditions, which
require that $B_{s, m}$ should be an integer for 
every $s \in \bbZ$; see \cite[Theorem 22.4.1]{hirzebruch}. 

\begin{lem} \label{lem.schw}
(a) $B_{s, m} = 0$ for every $s = 1, \dots, m - 2$. 

\smallskip
(b) $B_{1, k} = 0$ for $k = 3, \dots, m$.  
\end{lem}

\begin{proof}
(a) Since $B_{s, m}$ is an integer, it is enough to show 
that $|B_{s, m}| < 1$. Indeed, for each $i = 1, \dots, r$,
\begin{gather} 
\left| \begin{pmatrix} s - w_i \\ m \end{pmatrix} \right|
\le \frac{1}{m!} (|s - w_i| \dots |2 - w_i|) \cdot 
(|1 - w_i| \cdot |-w_i| \cdot |-1 - w_i|) \cdot 
\nonumber \\
\quad \quad \quad \quad \quad \quad \quad \quad \quad
(|-2 - w_i| \dots |s - m + 1 - w_i|)  \le \nonumber \\
\frac{1}{m!} \bigl( (s + 1) \cdot s \cdot \ldots \cdot3 \bigr) \cdot |1- w^2| 
\cdot \bigl( 3 \cdot 4 \cdot \ldots \cdot (m-s) \bigr) \le \nonumber \\ 
\frac{(s+1)! (m-s)!}{2m!} =   
\frac{(m + 1)}{2 \begin{pmatrix} m + 1 \\ s + 1 \end{pmatrix}} \le 
\frac{(m + 1)}{2 \begin{pmatrix} m + 1 \\ 2 \end{pmatrix}} = \frac{1}{m} 
\, . \nonumber \end{gather}
Note that we have used the inequality
$\begin{pmatrix} m + 1 \\ s + 1 \end{pmatrix} \ge 
\begin{pmatrix} m + 1 \\ 2 \end{pmatrix}$, which is valid for 
any $s = 1, \dots, m - 2$.  Now
\[ |B_{s, m}| \le \sum_{i = 1}^r \left|
\begin{pmatrix} s - w_i \\ m \end{pmatrix} \right| 
 \le r \frac{1}{m} < 1 \, , \]
and part (a) follows.

\smallskip
(b) Combining part (a) with the identity 
$B_{s, m} - B_{s-1, m} = B_{s-1, m-1}$, we conclude that
$B_{i, m-1} = 0$ for $i = 1, \dots, m - 3$. Repeating this argument, 
we see that for every $j = 0, \dots, m - 3$, and every
$i = 1, \dots, m - 2 - j$, we have $B_{i, m - j} = 0$. Now
set $i = 1$ and $k = m-j$, and part (b) follows.
\end{proof}

We now return to the proof of Proposition~\ref{prop.schw}.
By the  Chinese Remainder theorem, the semisimple 
$\bbQ$-algebra $F = \bbQ[x]/(p(x))$ is isomorphic
to $\bbQ(w_1) \oplus \dots \oplus \bbQ(w_d)$ via
$p(x) \mapsto (p(w_1), \dots, p(w_d))$. Thus for every
$f(x) \in F$, 
\[ \tr_{F/\bbQ} \, f(x) = \sum_{i = 1}^r f(w_i) \, . \]  
In particular, setting $a = - (1 - x) x (-1-x) \in F$,
and $b_0 = 1$, $b_1 = -2 - x$, $b_2 = (-2-x)(-3-x)$, $\dots$,
$b_{m-3} =  (-2 - x)(-3-x) \dots (2 - m - x)$ in $F$, 
and applying Lemma~\ref{lem.schw}(b), we see that
\begin{equation} \label{e.schw3}
\text{$\tr_{F/\bbQ}(a b_i) = (i+3)! B_{1, i + 3} = 0$ 
for every $i = 0, 1, \dots, m - 3$.} 
\end{equation}
Note since $1, x, \dots, x^{r-1}$ are $\bbQ$-linearly 
independent in $F$, so are $b_0, \dots, b_{r-1}$. 
Since we are assuming $r  \le m - 2$, 
this implies that $b_0, b_1, \dots, b_{m-3}$ span $F$. But the 
trace form $\mathopen< {\bf x}, {\bf y} \mathclose> =
\tr_{F/\bbQ}(xy)$ is non-singular on $F$ (viewed as an $r$-dimensional
$\bbQ$-vector space); thus~\eqref{e.schw3} is only possible if
$a = 0$ in $F$.  In other words, $a(w_i) = 
(1- w_i) (- w_i) (-1- w_i) = 0$ or, equivalently, 
$w_i = \pm 1$ for every $i = 1, \dots, r$. Since $w_1, \dots, w_r$
are assumed to be distinct, we conclude that either $r = 1$ and 
$w_1 = \pm 1$ or $r = 2$ and $\{w_1, w_2 \} = \{-1, 1 \}$. This completes
the proof of Proposition~\ref{prop.schw}
\end{proof}

\section{Proof of Theorem~\ref{thm1b}}

From now on we will assume that $i = 1$. That is, we will be
interested in nesting maps $\bbP_F^{n-1} \lra \Gr(j, F^n)$,
where $2 \le j \le n-1$. 

\begin{lem} \label{lem.1b} Suppose there exists a continuous nesting map
$f \colon \bbP_{\bbC}^{n-1} \lra \Gr(j, \bbC^n)$ (respectively,
an algebraic nesting map $f \colon \bbP_{F}^{n-1} \lra \Gr(j, F^n)$), where
$2 \le j \le n-1$ and $F$ is an algebraically closed field.   

\smallskip
(a) There exists topological
(respectively, algebraic) vector bundles $\cB$ and $\cC$ on
$\bbP_{\bbC}^{n-1}$ (respectively, $\bbP_F^{n-1}$) of ranks
$n - j$ and $j - 1$ with Chern polynomials $p(t)$ and $q(t) \in \bbZ[t]$
such that $p(t) q(t) = 1 + t + \dots + t^{n-1}$. In fact, we can take
$\cC = \cA/\cB$, where $\cA$ and $\cB$ are defined in~\eqref{e.AB}.

\smallskip
(b) $n$ is even and either $j = n-1$ and $p(t) = t + 1$ or
$j = 2$ and $q(t) = t + 1$.
\qed
\end{lem}   

\begin{proof} (a) We specialize the argument of 
Section~\ref{sect.i>=2} to the case where $i = 1$. In this case
the cohomology ring $H^{\ast}(\bbP^{n-1}, F^n) =
H^{\ast}(\Gr(i, F^n), \bbZ)$ reduces to $\bbZ[h]/(h^n = 0)$,
where $h$ is the class of a hyperplane section in $\bbP^{n-1}$. (In 
Section~\ref{sect.i>=2} we denoted $h$ by $\phi(\sigma_1)$.)
Defining $\cA$ and $\cB$ as in~\eqref{e.AB} (with $i = 1$),
setting $\cC = \cA/\cB$, and writing $t$ for $x_1$, we see that
\[ p(t) q(t) = \phi(c_{\rm tot}(\cA)) = 
1 + \theta_1(t) + \dots + \theta_{n-1}(t) 
= 1 + t + \dots + t^{n-1} \, ; \]
cf.~\eqref{e2.5}.

\smallskip
(b) If $n = 3$, the polynomial  
$1 + t + t^2$ is irreducible, contradicting part (a).
Thus we may assume $n \geq 4$. In this case, 
$1 \leq n -j \leq n-3$ or $1 \leq j - 1 \leq n-3$, and thus
Proposition~\ref{prop.schw} (with $m = n-1$)
applies to the bundle $\cB$ and to the bundle $\cC$ of part (a).
Since $p(t)$ and $q(t)$ have no multiple factors, their roots are
$n$th roots of unity, and $p(1), q(1) \neq 0$,
we see that the only possibilities for the Chern polynomials
are the ones listed in part (b).
\end{proof}

We now turn to the proof of Theorem~\ref{thm1b}. Part (a)
is an immediate consequence of the above lemma. To prove part (b)
of Theorem~\ref{thm1b}, assume that there exists an algebraic nesting map 
\[ f \colon \bbP_{F}^{n-1} \lra \Gr(j, F^n) \, . \]
By Lemma~\ref{lem.1b}(b), $n$ is even and $j = 2$ or $n - 1$.  
Our goal is to show that (i) $j = 2$ is impossible, and
(ii) if $j = n-1$ then $f(l) = l^{\perp_{\omega}}$ for some alternating
form $\omega$ on $F^n$, as in Example~\ref{ex1}.
(For $F = \bbC$, (ii) was proved in by Roan~\cite{roan}; we will give 
a short characteristic-free proof below.) 

\smallskip
Proof of (i): Suppose $j = 2$. Consider the exact sequences
\begin{gather} \label{e5.1a}
0 \lra \cO_{\bbP^{n-1}}(-1) \lra \cO_{\bbP^{n-1}}^{\oplus n} \lra \cA \lra 0 
\nonumber \\
\label{e5.1b}
0 \lra \cB \lra \cA \lra \cC \lra 0 \nonumber 
\end{gather}
of algebraic vector bundles on $\bbP_F^{n-1}$.
From the first sequence, we see that $H^0(\bbP_F^{n-1}, \cA) = F^n$.
By Lemma~\ref{lem.1b}(b), $\cB$ is a line bundle with
Chern polynomial $1 + t$, i.e., $\cB = \cO(1)$ on $\bbP_F^{n-1}$.
Thus $H^0(\bbP_F^{n-1}, \cB) = F^n$, and the second sequence
yields
\[ (0) \lra F^n \lra F^n \lra H^0(\bbP_F^{n-1}, \cC) \lra (0) \, ,  \]
which means that $\cC$ has no global sections. On the other hand,
$\cC$, being a quotient of $\cA$, is generated by global sections.
This contradiction shows that $f$ cannot exist for $j = 2$.

\smallskip 
Proof of (ii): Assume $j = n-1$. Consider the exact sequence
\[ 0 \lra \cT_{n-1} \lra \cO_{\check{\bbP}^{n-1}}^{\oplus n} 
\lra L \lra 0 \, , \]
where $\cT_{n-1}$ is the tautological bundle on 
$\Gr(n-1, F^n) = \check{\bbP}^{n-1}$.  Pulling back this sequence to 
$\bbP^{n-1}$ via $f$, we obtain an exact sequence 
\[ 0 \lra f^{\ast}(\cT_{n-1}) \lra \cO_{\bbP^{n-1}}^{\oplus n}
\lra f^{\ast}(L) \lra 0 \]
of vector bundles on $\bbP^{n-1}$.
Since $L$ is a line bundle generated by global sections, 
$L = \cO_{\check{\bbP}^{n-1}}(d)$ for some $d > 0$. 
Since $f$ is nesting,  
$\cT_1 = \cO_{\bbP^{n-1}}(-1) \subset f^{\ast}(\cT_{n-1})$. Recall that
the bundles $\cB$ and $\cC$ in 
the statement of Lemma~\ref{lem.1b} are defined by
$\cB = f^{\ast}(\cT_{n-1})/\cT_1$ and $\cC = \cA/\cB$, where
$\cA = \cO_{\bbP^{n-1}}^{\oplus n}/\cT_1$; cf.~\eqref{e.AB}.
Thus the line bundle $\cC$ is isomorphic 
to $\cO_{\bbP^{n-1}}/f^{\ast}(\cT_{n-1}) = f^{\ast}(L)$. 
By Lemma~\ref{lem.1b}(b), the Chern polynomial of
$\cC$ is $1 + t$; in other words, 
$\cC = \cO_{\bbP^{n-1}}(1)$. On the other hand,
$\cC = f^{\ast}(L) = f^{\ast}(\cO_{\check{\bbP}^{n-1}}(d))$. This 
is only possible if $d = 1$ and $f$ is induced by a non-singular linear map 
\[ f^{\ast} \colon H^0(\check{\bbP}^{n-1}, \cO_{\check{\bbP}^{n-1}}(1)) \lra
H^0(\bbP^{n-1}, \cO_{\bbP^{n-1}}(1)) \,  \]
In other words, $f$ is induced by a non-singular linear map
$f^{\ast} \colon F^n \lra (F^n)^{\ast}$, where $(F^n)^*$ 
is the dual vector space to $F^n$, $\bbP^{n-1} = \bbP(F^n)$ and 
$\check{\bbP}^{n-1} = \bbP((F^n)^*)$.
Now $\omega(x, y) = f^{\ast}(x)(y)$ 
is a non-singular bilinear form on $F^n$. Since $f$ is nesting,
$\omega(x, x) = 0$ for every $x \in F^n$, i.e., $\omega$ is alternating.
This shows that $f$ is obtained by the construction 
in Example~\ref{ex1}, as claimed. The proof of 
Theorem~\ref{thm1b} is now complete.
\qed

\section{Proof of Corollary~\ref{cor2}}

Corollary~\ref{cor2} is an immediate consequence of Theorem~\ref{thm1b}
and the lemma below.

\begin{lem} \label{lem4.1} The following are in natural
(i.e., $\PGL_n$-equivariant) correspondence:

\smallskip
(a) algebraic nesting maps
$f \colon \bbP_{F}^{n-1} \lra \Gr(j, F^n)$,

\smallskip
(b) algebraic subbundles of rank $j - 1$ of
the bundle $\cA = \cO^{n+1}/\cO(-1)$ on $\bbP_{F}^{n-1}$,

\smallskip
(c) algebraic subbundles of rank $j - 1$ of the tangent bundle
$T({\bbP_{F}^{n-1}})$,

\smallskip
\noindent
For $F = \bbR$ or $\bbC$ the lemma remains true if ``algebraic"
is replaced by ``continuous". 
\end{lem}

\begin{proof} 
Let $\cT_j$ be the tautological bundle on
$\Gr(j, F^n)$ (as in Section~\ref{sect.i>=2}); in
particular, $\cT_1 = \cO(-1)$ is the tautological line bundle
on $\bbP_{F}^{n-1}$.  Given a nesting map
$f \colon \bbP_{F}^{n-1} \lra \Gr(j, F^n)$, we associate to it the
subbundle $\cB = f^{\ast}(\cT_j)/\cO(-1)$ of $\cA$ of rank $j - 1$;
cf.~\eqref{e.AB}. 

Conversely, a subbundle $\cB$ of $\cA$ of rank $j -1$ lifts to
a subbundle $\cB'$ of $\cO^{n}$ of rank $j$ containing $\cO(-1)$, 
and we can define a nesting map
\[ f \colon \bbP_{F}^{n-1} \lra \Gr(i, F^n) \]
by $p \mapsto \cB'(p)$.
This establishes a natural bijective correspondence between (a) and (b).

To show that (b) and (c) are in a natural bijective
correspondence, note that $\cA = T({\bbP_F^{n-1}})(-1)$;
see, e.g.,~\cite[p. 6]{oss} or \cite[p. 409]{gh}. The same argument
works in the continuous case.
\end{proof}

\section{Appendix: Proof of Theorem~\ref{thm0}}
\label{sect2}

For every $l \in \Gr(i, F^n)$, let
\[ X_l = \{ L \in \Gr(j, F^n) \ | \, l \subset L \} \, . \]
Since the sets $\Gr(i, F^n)$ and $\Gr(j, F^n)$ have
the same cardinality, a bijective nesting map
$f \colon \Gr(i, F^n) \lra \Gr(j, F^n)$ may be viewed as
a system of distinct representatives for the collection of subsets
$\{ X_l \, | \, l \in \Gr(i, F^n) \}$ of $\Gr(j, F^n)$. Indeed, given
a system of distinct representatives $\{ x_l \}$,
with $x_l \in X_l$, the map defined by $f(l) = x_l$ is nesting and bijective.
Conversely, given $f$, as in Theorem~\ref{thm0}, the elements $x_l = f(l)$
form a system of distinct representatives for $\{ X_l \}$.

Thus, by a theorem of P. Hall (see, e.g.,~\cite[Theorem 5.1.1]{ryser}),
we only need to check that
\begin{equation} \label{e.hall}
|X_{l_1} \cup  \dots \cup  X_{l_k}| \geq k
\end{equation}
for every choice of distinct elements $l_1, \dots, l_k \in \Gr(i, F^n)$.

Let $N$ be the number of $j$-planes in $F^n$ containing a given $i$-plane.
Since $i + j = n$, 
\[ N = | \Gr(j - i, F^{n-i}) | = | \Gr(n-j, F^{n-i}) | =
| \Gr(i, F^{j}) | \]
is also the number of $i$-planes in $F^n$ contained in a given $j$-plane.
To prove~\eqref{e.hall}, we will count the number of elements in the set
\[ \text{$W = \{ (l, L) \, |$ where $l \subset L$ and $l = l_1, \dots, l_k \,
  \}$} \]
in two ways. On the one hand, projecting to the first component, we see
that $|W| = k N$. On the other hand, projecting to the second component,
we obtain $|W| \leq |X_{l_1} \cup  \dots \cup  X_{l_k}| N$. Thus
\[ |W| = k N \leq
|X_{l_1} \cup  \dots \cup  X_{l_k}| N \]
and~\eqref{e.hall} follows.
\qed

\bibliographystyle{amsplain}

\providecommand{\bysame}{\leavevmode\hbox to3em{\hrulefill}\thinspace}

\end{document}